\documentclass[fleqn]{mat01}
\usepackage{times,mathtimy,amssymb,epsfig}
\begin{document}

\setcounter{page}{379}
\firstpage{379}

\def\thoe{\trivlist\item[\hskip\labelsep{{\bf Theorem}}]}
\newtheorem{theo}{\bf Theorem}
\renewcommand\thetheo{\arabic{theo}}
\newtheorem{theor}[theo]{\bf Theorem}
\newtheorem{propo}[theo]{\rm PROPOSITION}
\newtheorem{coro}[theo]{\rm COROLLARY}
\newtheorem{lem}[theo]{Lemma}
\newtheorem{fact}[theo]{Fact}
\newtheorem{claim}{Claim}
\newtheorem{rem}[theo]{Remark}
\newtheorem{definit}[theo]{\rm DEFINITION}
\newtheorem{exampl}{Example}
\def\maintheorem{\trivlist\item[\hskip\labelsep{{\bf Main Theorem.}}]}
\def\pmt{\trivlist\item[\hskip\labelsep{{\it Proof of Main Theorem.}}]}

\markboth{Liu Lanzhe}{Triebel--Lizorkin space estimates}

\title{Triebel--Lizorkin space estimates
for multilinear operators\\ of sublinear operators}

\author{LIU LANZHE}

\address{Department of Applied Mathematics, Hunan University,
Changsha 410082,\\
\noindent People's Republic of China\\
\noindent Email: lanzheliu@263.net}

\volume{113}

\mon{November}

\parts{4}

\Date{MS received 23 March 2003}

\begin{abstract}
In this paper, we obtain the continuity for some multilinear operators
related to certain non-convolution operators on the Triebel--Lizorkin
space. The operators include Littlewood--Paley operator and Marcinkiewicz
operator.
\end{abstract}


\maketitle

\section{Introduction}

Let $T$ be the singular integral operator, a well-known result of
Coifman, Rochberg and Weiss \cite{6} which states that the commutator
$[b,T]=T(bf)-bTf$ (where $b\in \,\hbox{BMO}$) is bounded on
$L^p(R^n) (1<p<\infty)$. Chanillo \cite{1} proves a similar result
when $T$ is replaced by the fractional integral operator. In \cite{9,11},
these results on the Triebel--Lizorkin spaces and the case $b\in
\hbox{Lip}\beta$ (where $\hbox{Lip}\beta$ is the homogeneous Lipschitz space) are
obtained. The main purpose of this paper is to study the continuity for
some multilinear operators related to certain convolution operators on
the Triebel--Lizorkin spaces. In fact, we shall obtain the continuity on
the Triebel--Lizorkin spaces for the multilinear operators only under
certain conditions on the size of the operators. As applications, we
prove the continuity of the multilinear operators related to the
Littlewood--Paley operator and Marcinkiewicz operator on the
Triebel--Lizorkin spaces.

\section{Notations and results}

In the sequel, $Q$ will denote a cube of $R^n$ with sides parallel to
the axes, and for a cube $Q$, let $f_Q=|Q|^{-1}\int_Qf(x)\hbox{d}x$ and
$f^{\#}(x)=\sup_{x\in Q}|Q|^{-1}\int_Q|f(y)-f_Q|\hbox{d}y$. For
$1\le r<\infty$ and $0\le \delta$, let
\begin{equation*}
M_{\delta,r}(f)(x)=\sup_{x\in Q}\left(\frac{1}{|Q|^{1-\delta
r/n}}\int_Q|f(y)|^r\hbox{d}y\right)^{1/r}.
\end{equation*}
We denote $M_{\delta,r}(f)=M_r(f)$ if $\delta=0$, which is the
Hardy--Littlewood maximal function when $r=1$ (in this case, we denote
$M_1(f)=M(f)$, see \cite{12,13}). For $\beta>0$ and $p>1$, let $\dot
F_p^{\beta,\infty}$ be the homogeneous Triebel--Lizorkin space, the
Lipschitz space $\dot\land_\beta$ is the space of functions $f$ such
that
\begin{equation*}
\Vert f\Vert_{\dot\land_\beta}=\sup_{\substack{x,h\in R^n \\
h\ne 0}} |\Delta_h^{[\beta]+1}f(x)|/|h|^\beta<\infty,
\end{equation*}
where $\Delta_h^k$ denotes the $k$th difference operator \cite{11}.

We are going to consider the multilinear operator as follows:

Let $m$ be a positive integer and $A$ a function on $R^n$. We denote
\begin{equation*}
R_{m+1}(A;x,y)=A(x)-\sum_{|\alpha|\le m}\frac{1}{\alpha!}D^\alpha
A(y)(x-y)^\alpha.
\end{equation*}

\begin{definit}$\left.\right.$\vspace{.5pc}

\noindent {\rm Define $F(x,y,t)$ on $R^n\times R^n\times [0,+\infty)$.
Then we denote
\begin{equation*}
F_t(f)(x)=\int_{R^n}F(x,y,t)f(y)\hbox{d}y
\end{equation*}
and
\begin{equation*}
F_t^A(f)(x)=\int_{R^n}\frac{R_{m+1}(A;x,y)}{|x-y|^m}F(x,y,t)f(y)\hbox{d}y.
\end{equation*}
Let $H$ be the Hilbert space $H=\{h:\Vert h\Vert <\infty\}$. For each fixed
$x\in R^n$, we view $F_t(f)(x)$ and $F_t^A(f)(x)$ as a mapping from
$[0,+\infty)$ to $H$. Then, the multilinear operators related to $F_t$
is defined by
\begin{equation*}
T^Af(x)=\Vert F_t^A(f)(x)\Vert ;
\end{equation*}
we also define that $Tf(x)=\Vert F_tf(x)\Vert $.}
\end{definit}

In particular, we consider the following two sublinear operators. Fix
$\lambda > 1$.

\begin{definit}$\left.\right.$\vspace{.5pc}

\noindent {\rm Let $\varepsilon>0$ and $\psi$ be a
  fixed function which satisfies the following properties:
\begin{enumerate}
\renewcommand{\labelenumi}{(\arabic{enumi})}
\item $|\psi(x)|\le C(1+|x|)^{-(n+1)}$,

\item $|\psi(x+y)-\psi(x)|\le
C|y|^\varepsilon(1+|x|)^{-(n+1+\varepsilon)}$ when $2|y|<|x|$.
\end{enumerate}

The multilinear Littlewood--Paley operator is defined by
\begin{equation*}
g_\lambda^A(f)(x)=\left[\int\int_{R_+^{n+1}}
\left(\frac{t}{t+|x-y|}\right)^{n\lambda} |F_t^A(f)(x,y)|^2
\frac{{\rm d}y{\rm d}t}{t^{n+1}}\right]^{1/2},
\end{equation*}
where
\begin{equation*}
F_t^A(f)(x,y)=\int_{R^n}\frac{R_{m+1}(A;x,z)}{|x-z|^m}
f(z)\psi_t(y-z){\rm d}z
\end{equation*}
and $\psi_t(x)=t^{-n}\psi(x/t)$ for $t>0$. We denote
$F_t(f)(y)=f\ast\psi_t(y)$. We also define
\begin{equation*}
g_\lambda(f)(x)=\left(\int\int_{R_+^{n+1}}
\left(\frac{t}{t+|x-y|}\right)^{n\lambda} |F_t(f)(y)|^2
\frac{{\rm d}y{\rm d}t}{t^{n+1}}\right)^{1/2},
\end{equation*}
which is the Littlewood--Paley operator \cite{13}.

Let $H$ be the Hilbert space $H= \{h:\Vert h\Vert =
(\int\int_{R_+^{n+1}} |h(t)|^2{\rm d}y{\rm d}t/t^{n+1})^{1/2}<\infty\}$.
Then for each fixed $x\in R^n$, $F_t^A(f)(x,y)$ may be viewed as a
mapping from $(0, +\infty)$ to $H$, and it is clear that
\begin{align*}
g_\lambda^A(f)(x) &= \left\Vert\left(\frac{t}{t+|x-
y|}\right)^{n\lambda/2} F_t^A(f)(x,y)\right\Vert,\\[.2pc]
g_\lambda(f)(x) &= \left\Vert\left(\frac{t}{t+|x-y|}\right)^{n\lambda/2}
F_t(f)(y) \right\Vert.
\end{align*}}
\end{definit}\vspace{.5pc}

\begin{definit}$\left.\right.$\vspace{.5pc}

\noindent {\rm Let $0<\gamma\le 1$ and $\Omega$ be a homogeneous of
degree zero on $R^n$ such that $\int_{S^{n-1}}\Omega(x'){\rm
d}\sigma(x')=0$. Assume that $\Omega\in {\rm Lip}_\gamma(S^{n-1})$, that
is there exists a constant $M>0$ such that for any $x,y\in S^{n-1}$,
$|\Omega(x)-\Omega(y)|\le M|x-y|^\gamma$. We denote
$\Gamma(x)=\{(y,t)\in R_+^{n+1}:|x-y|<t\}$ and the characteristic
function of $\Gamma(x)$ by $\chi_{\Gamma(x)}$. The multilinear
Marcinkiewicz integral operator is defined by
\begin{equation*}
\mu_\lambda^A (f)(x)=\left[\int\int_{R_+^{n+1}}
\left(\frac{t}{t+|x-y|}\right)^{n\lambda} |F_t^A(f)(x)|^2
\frac{{\rm d}y{\rm d}t}{t^{n+3}}\right]^{1/2},
\end{equation*}
where
\begin{equation*}
F_t^A(f)(x,y)=\int_{|y-z|\le t}\frac{\Omega(y-z)}{|y-z|^{n-1}}
\frac{R_{m+1}(A;x,z)}{|x-z|^m}f(z){\rm d}z.
\end{equation*}
We denote that
\begin{equation*}
F_t(f)(y)=\int_{|y-z|\le t}\frac{\Omega(y-z)}{|y-z|^{n-1}}
f(z){\rm d}z.
\end{equation*}
We also define that
\begin{equation*}
\mu_\lambda(f)(x)=\left(\int\int_{R_+^{n+1}}\left(\frac{t}{t+|x-
y|}\right)^{n\lambda} |F_t(f)(y)|^2 \frac{{\rm d}y{\rm
d}t}{t^{n+3}}\right)^{1/2},
\end{equation*}
which is the Marcinkiewicz integral operator \cite{14}.

Let $H$ be the Hilbert space $H= \{h: \Vert h\Vert = (
\int\int_{R_+^{n+1}}|h(t)|^2{\rm d}y{\rm d}t/t^{n+3})^{1/2}
<\infty\}$. Then for each fixed $x\in R^n$, $F_t^A(f)(x,y)$ may be
viewed as a mapping from $(0, +\infty)$ to $H$, and it is clear that
\begin{align*}
\mu_\lambda^A(f)(x) &= \left\Vert\left(\frac{t}{t+|x-
y|}\right)^{n\lambda/2} F_t^A(f)(x,y)\right\Vert,\\[.2pc]
\mu_\lambda(f)(x) &= \left\Vert\left(\frac{t}{t+|x-
y|}\right)^{n\lambda/2} F_t(f)(y) \right\Vert.
\end{align*}
}
\end{definit}

It is clear that Definitions~2 and 3 are the particular examples of
Definition 1. Note that when $m=0, T_A$ is just the commutator of
$F_t$ and $A$, and when $m>0$, it is the non-trivial generalizations of
the commutators. It is well-known that multilinear operators are of
great interest in harmonic analysis and have been widely studied by many
authors (see \cite{2,3,4,5,7,8}). The main purpose of this paper is to study
the continuity for the multilinear operators on the Triebel--Lizorkin
spaces. We shall prove the following theorems in \S3.

\setcounter{theo}{0}
\begin{theor}[\!] Let $g_\lambda^A$ be the multilinear Littlewood--Paley
operator as in Definition~{\rm 2}. If $0<\beta<1/2$ and $D^\alpha
A\in\dot\land_\beta$ for $|\alpha|=m$. Then
\begin{enumerate}
\renewcommand{\labelenumi}{{\rm (\alph{enumi})}}
\item $g_\lambda^A$ maps $L^p(R^n)$ continuously onto $\dot
F_p^{\beta,\infty}(R^n)$ for $1<p<\infty${\rm ;}

\item $g_\lambda^A$ maps $L^p(R^n)$ continuously onto $L^q(R^n)$ for
$1<p<n/\beta$ and $1/p-1/q=\beta/n$.
\end{enumerate}
\end{theor}

\begin{theor}[\!] Let $\mu_\lambda^A$ be the multilinear Marcinkiewicz
operator as in Definition~{\rm 3}. If $0<\beta<1/2$ and $D^\alpha
A\in\dot\land_\beta$ for $|\alpha|=m$. Then
\begin{enumerate}
\renewcommand{\labelenumi}{{\rm (\alph{enumi})}}
\item $\mu_\lambda^A$ maps $L^p(R^n)$ continuously onto $\dot
F_p^{\beta,\infty}(R^n)$ for $1<p<\infty${\rm ;}

\item $\mu_\lambda^A$ maps $L^p(R^n)$ continuously onto $L^q(R^n)$ for
$1<p<n/\beta$ and $1/p-1/q=\beta/n$.
\end{enumerate}\vspace{-2pc}
\end{theor}

\section{Main theorem and proof}

We first prove a general theorem.

\begin{maintheorem} {\it Let $0<\beta<1$ and $D^\alpha A\in\dot\land_\beta$
for $|\alpha|=m$. Suppose $F_t, T$ and $T^A$ are the same as in
Definition~{\rm 1,} if $T$ is bounded on $L^p(R^n)$ for $1<r<\infty$ and $T$
satisfies the following size condition{\rm :}
\begin{equation*}
\Vert F_t^A(f)(x)-F_t^A(f)(x_0)\Vert \le C\sum_{|\alpha|=m}\Vert
D^\alpha A\Vert_{\dot\land_\beta}|Q|^{\beta/n}M(f)(x)
\end{equation*}
for any cube $Q$ with ${\rm supp}\, f \subset(2Q)^c$ and $x\in Q$. Then

\begin{enumerate}
\renewcommand{\labelenumi}{{\rm (\alph{enumi})}}
\item $T^A$ maps $L^p(R^n)$ to $\dot F_p^{\beta, \infty}(R^n)$ for
$1<p<\infty${\rm ;}

\item $T^A$ maps $L^p(R^n)$ to $L^q(R^n)$ for $1<p<n/\beta$ and
$1/q=1/p-\beta/n$.
\end{enumerate}}
\end{maintheorem}

To prove the theorem, we need the following lemmas.

\setcounter{theo}{0}
\begin{lem}\hskip -.5pc{\rm \cite{11}}.\ \ For $0<\beta<1, 1<p<\infty${\rm ,} we have
\begin{align*}
\Vert f\Vert_{\dot F_p^{\beta,\infty}} &\approx \left\Vert\sup_Q
\frac{1}{|Q|^{1+\beta/n}} \int_Q|f(x)-f_Q|{\rm d}x\right\Vert_{L^p}\\[.2pc]
&\approx \left\Vert\sup_{\cdot\in Q}
\inf_c\frac{1}{|Q|^{1+\beta/n}}\int_Q|f(x)-c|{\rm d}x\right\Vert_{L^p}.
\end{align*}
\end{lem}

\pagebreak

\begin{lem}\hskip -.5pc{\rm \cite{11}}.\ \ For $0<\beta<1, 1\le p\le \infty${\rm ,} we have
\begin{align*}
\Vert f\Vert_{\dot \land_\beta} &\approx
\sup_Q\frac{1}{|Q|^{1+\beta/n}}\int_Q |f(x)-f_Q|{\rm d}x\\[.2pc]
&\approx \sup_Q \frac{1}{|Q|^{\beta/n}}\left(\frac{1}{|Q|}\int_Q|f(x)-
f_Q|^p{\rm d}x\right)^{1/p}.
\end{align*}
\end{lem}

\begin{lem}\hskip -.5pc{\rm \cite{1,2}}.\ \ Suppose that $1\le r<p<n/\delta$ and
$1/q=1/p-\delta/n$. Then
\begin{equation*}
\Vert M_{\delta,r}(f)\Vert_{L^q}\le C\Vert f\Vert_{L^p}.
\end{equation*}
\end{lem}

\begin{lem}\hskip -.5pc{\rm \cite{5}}.\ \ Let $A$ be a function on $R^n$ and $D^\alpha A\in
L^q(R^n)$ for $|\alpha|=m$ and some $q>n$. Then
\begin{equation*}
|R_m(A ; x,y)|\le C|x-y|^m\sum_{|\alpha|=m}\left(\frac{1}{|\tilde Q
(x,y)|}\int_{\tilde Q(x,y)}|D^\alpha A(z)|^q {\rm d}z\right)^{1/q},
\end{equation*}

$\left.\right.$\vspace{-1.5pc}

\noindent where $\tilde Q(x,y)$ is the cube centered at $x$ and having side length
$5\sqrt{n}|x-y|$.\vspace{.5pc}
\end{lem}

\begin{pmt}$\left.\right.$

\noindent (a) Fix a cube $Q=Q(x_0,l)$ and $\tilde x\in Q$. Let $\tilde
Q=5\sqrt{n}Q$ and $\tilde
A(x)=A(x)-\sum_{|\alpha|=m}\frac{1}{\alpha!} (D^\alpha A)_{\tilde
Q}x^\alpha$, then $R_m(A;x,y)=R_m(\tilde A;x,y)$ and $D^\alpha\tilde
A=D^\alpha A -(D^\alpha A)_{\tilde Q}$ for $|\alpha|=m$. We write, for
$f_1=f\chi_{\tilde Q}$ and $f_2=f\chi_{R^n\setminus\tilde Q}$,
\begin{align*}
F_t^A(f)(x) &= \int_{R^n}\frac{R_{m+1}(\tilde
A;x,y)}{|x-y|^m}F(x,y,t)f(y){\rm d}y\\[.2pc]
&= \int_{R^n}\frac{R_{m+1}(\tilde A;x,y)}{|x-y|^m}F(x,y,t)f_2(y){\rm d}y\\[.2pc]
&\quad\ +\int_{R^n}\frac{R_m(\tilde A;x,y)}{|x-y|^m}F(x,y,t)f_1(y){\rm d}y\\[.2pc]
&\quad\ -\sum_{|\alpha|=m}\frac{1}{\alpha!}
\int_{R^n}\frac{F(x,y,t)(x-y)^\alpha}{|x-y|^m}D^\alpha\tilde
A(y)f_1(y){\rm d}y,
\end{align*}
then
\begin{align*}
&|T^A(f)(x)-T^{\tilde A}(f_2)(x_0) | = | \Vert
F_t^A(f)(x)\Vert -\Vert F_t^{\tilde A}(f_2)(x_0)\Vert|\\[.2pc]
&\le \left\Vert F_t\left(\frac{R_m(\tilde
A;x,\cdot)}{|x-\cdot|^m}f_1\right)(x)\right\Vert
+\sum_{|\alpha|=m}\frac{1}{\alpha!}\left\Vert F_t\left(\frac{(x-
\cdot)^\alpha}{|x-\cdot|^m}D^\alpha\tilde A f_1\right)(x)\right\Vert\\[.2pc]
&\quad\ + \Vert F_t^{\tilde A}(f_2)(x)-F_t^{\tilde A}(f_2)(x_0)\Vert
={\rm I}(x) + {\rm II}(x) + {\rm III}(x).
\end{align*}
Thus,
\begin{align*}
&\frac{1}{|Q|^{1+\beta/n}}\int_Q\left|T^A(f)(x)- T^{\tilde
A}(f)(x_0)\right|{\rm d}x \\[.2pc]
&\le \frac{1}{|Q|^{1+\beta/n}}\int_Q {\rm I}(x){\rm
d}x+\frac{1}{|Q|^{1+\beta/n}} \int_Q {\rm II}(x){\rm
d}x+\frac{1}{|Q|^{1+\beta/n}}\int_Q {\rm III}(x){\rm d}x\\[.2pc]
&:= {\rm I} + {\rm II} + {\rm III}.
\end{align*}
Now, let us estimate I, II and III, respectively. First, for $x\in
Q$ and $y\in \tilde Q$, using Lemmas~2 and 4, we get
\begin{align*}
|R_m(\tilde A; x,y)| &\le C|x-y|^m\sum_{|\alpha|=m}\sup_{x\in \tilde
Q}|D^\alpha A(x)-(D^\alpha A)_{\tilde Q}|\\[.2pc]
&\le C|x-y|^m|Q|^{\beta/n}\sum_{|\alpha|=m}\Vert D^\alpha
A\Vert_{\dot\land_\beta}.
\end{align*}
Thus, by Holder's inequality and the $L^r$ boundedness of $T$ for
$1<r<p$, we obtain
\begin{align*}
{\rm I} &\le C\sum_{|\alpha|=m}\Vert D^\alpha
A\Vert_{\dot\land_\beta}\frac{1}{|Q|}\int_Q|T(f_1)(x)|{\rm d}x\\[.2pc]
&\le C\sum_{|\alpha|=m}\Vert D^\alpha A\Vert_{\dot\land_\beta}\Vert
T(f_1)\Vert_{L^r}|Q|^{-1/r} \\[.2pc]
&\le C\sum_{|\alpha|=m}\Vert D^\alpha A\Vert_{\dot\land_\beta}\Vert
f_1\Vert_{L^r}|Q|^{-1/r}\\[.2pc]
&\le C\sum_{|\alpha|=m}\Vert D^\alpha
A\Vert_{\dot\land_\beta}M_r(f)(\tilde x).
\end{align*}
Secondly, for $1<r<q$, using the following inequality \cite{11}:
\begin{equation*}
\Vert D^\alpha A-(D^\alpha A)_{\tilde Q}f\chi_{\tilde Q}\Vert_{L^r}\le
C|Q|^{1/r+\beta/n} \Vert D^\alpha A\Vert_{\dot\land_\beta}M_r(f)(x)
\end{equation*}
and similar to the proof of I, we gain
\begin{align*}
{\rm II} &\le \frac{C}{|Q|^{1+\beta/n}}\sum_{|\alpha|=m}\Vert T((D^\alpha
A-(D^\alpha A)_{\tilde Q})f\chi_{\tilde Q})\Vert_{L^r}|Q|^{1-1/r} \\[.2pc]
&\le C|Q|^{-\beta/n-1/r}\sum_{|\alpha|=m}\Vert (D^\alpha A-(D^\alpha
A)_{\tilde Q})f\chi_{\tilde Q}\Vert_{L^r} \\[.2pc]
&\le C\sum_{|\alpha|=m}\Vert D^\alpha
A\Vert_{\dot\land_\beta}M_r(f)(\tilde x).
\end{align*}
For III, using the size condition of $T^A$, we have
\begin{equation*}
{\rm III}\le C\sum_{|\alpha|=m}\Vert D^\alpha A\Vert_{\dot\land_\beta}
M(f)(\tilde x).
\end{equation*}
We now put these estimates together, and taking the supremum over all
$Q$ such that $\tilde x\in Q$, and using the $L^p$ boundedness of $M_r$
for $r<p$, we obtain
\begin{equation*}
\Vert T^A(f)\Vert_{\dot F_p^{\beta,\infty}}\le C\sum_{|\alpha|=m}\Vert
D^\alpha A\Vert_{\dot\land_\beta}\Vert f\Vert_{L^p}.
\end{equation*}
This completes the proof of (a).

(b) By same argument as in the proof of (a), we have
\begin{align*}
&\frac{1}{|Q|}\int_Q |T^A(f)(x)-T^{\tilde A}(f_2)(x_0)|{\rm d}x\\[.2pc]
&\quad\ \le C\sum_{|\alpha|=m}\Vert D^\alpha A\Vert_{\dot\land_\beta}
(M_{\beta, r}(f)+ M_{\beta,1}(f)).
\end{align*}
Thus,
\begin{equation*}
(T^A(f))^{\#}\le C\sum_{|\alpha|=m}\Vert D^\alpha
A\Vert_{\dot\land_\beta}(M_{\beta,r}(f)+ M_{\beta,1}(f)).
\end{equation*}
Using Lemma~3, we gain
\begin{align*}
&\Vert T^A(f)\Vert_{L^q}\le C\Vert (T^A(f))^{\#}\Vert_{L^q}\\[.2pc]
&\quad\ \le C\sum_{|\alpha|=m}\Vert D^\alpha A\Vert_{\dot\land_\beta} (\Vert
M_{\beta, r}(f)\Vert_{L^q}+\Vert M_{\beta,1}(f)\Vert_{L^q})\le C\Vert
f\Vert_{L^p}.
\end{align*}
This completes the proof of (b) and the Main Theorem.
\end{pmt}

To prove Theorems~1 and 2, we need the following lemma.

\begin{lem} Let $1<p<\infty${\rm ,} $0<\beta<1$ and $D^\alpha
A\in\dot\land_\beta$ for $|\alpha|=m$. Then $g_\lambda^A$ and
$\mu_\lambda^A$ are all bounded on $L^p(R^n)$.
\end{lem}

\begin{proof} For $g_\lambda^A$, by Minkowski inequality and the
condition of $\psi$, we have
\begin{align*}
g_\lambda^A (f)(x) &\le \int_{R^n}\frac{|f(z)\Vert
R_{m+1}(A;x,z)|}{|x-z|^m}\\[.2pc]
&\quad\ \times \left(\int_{R_+^{n+1}}\left(\frac{t}{t+|x-y|}\right)^{n\lambda}
|\psi_t(y-z)|^2\frac{{\rm d}y{\rm d}t}{t^{1+n}}\right)^{1/2}{\rm d}z \\[.2pc]
&\le C\int_{R^n}\frac{|f(z)\Vert R_{m+1}(A;x,z)|}{|x-z|^m}\\[.2pc]
&\quad\ \times\left(\int_0^\infty\int_{R^n}\left(\frac{t}{t+|x-
y|}\right)^{n\lambda}
\frac{t^{-2n}}{(1+|y-z|/t)^{2n+2}} \right.\\[.2pc]
&\quad\ \left. \times \left(\frac{t}{t+|x-
y|}\right)^{n\lambda} \frac{{\rm d}y{\rm d}t}{t^{1+n}}\right)^{1/2}{\rm
d}z\\[.2pc]
&\le C\int_{R^n}\frac{|f(z)| |R_{m+1}(A;x,z)|}{|x-z|^m}\\[.2pc]
&\quad\ \times \left[\int_0^\infty\left(t^{-n}\!\int_{R^n}
\left(\frac{t}{t\!+\!|x\!-\!y|}\right)^{n\lambda}\!\! \frac{{\rm
d}y}{(t\!+\!|y\!-\!z|)^{2n+2}}\right)t{\rm d}t\right]^{1/2}{\rm d}z.
\end{align*}
Noting that
\begin{align*}
&t^{-n}\int_{R^n}\left(\frac{t}{t+|x-y|}\right)^{n\lambda}\frac{{\rm
d}y}{(t+|y- z|)^{2n+2}}\\[.2pc]
&\quad\ \le CM\left(\frac{1}{(t+|x-z|)^{2n+2}}\right) \le
\frac{C}{(t+|x-z|)^{2n+2}}
\end{align*}
and
\begin{equation*}
\int_0^\infty\frac{t{\rm d}t}{(t+|x-z|)^{2n+2}}=C|x-z|^{-2n},
\end{equation*}
we obtain
\begin{align*}
g_\lambda^A(f)(x) &\le C\int_{R^n}\frac{|f(z)\Vert
R_{m+1}(A;x,z)|}{|x-z|^m} \left( \int_0^\infty\frac{t{\rm
d}t}{(t+|x-z|)^{2n+2}}\right)^{1/2}{\rm d}z \\[.3pc]
&= C\int_{R^n}\frac{|f(z)\Vert R_{m+1}(A;x,z)|}{|x-z|^{m+n}}{\rm d}z.
\end{align*}
For $\mu_\lambda$, note that $|x-z|\le 2t$, $|y-z|\ge |x-z|-t\ge
|x-z|-3t$ when $|x-y|\le t$, $|y-z|\le t$, and $|x-z|\le t(1+2^{k+1})\le
2^{k+2}t$, $|y-z|\ge |x-z|-2^{k+3}t$ when $|x-y|\le 2^{k+1}t$, $|y-z|\le
t$, we have
\begin{align*}
\hskip -4pc\mu_\lambda^A (f)(x) &\le
\int_{R^n}\left[\int\int_{R_+^{n+1}}\left(\frac{t}
{t+|x-y|}\right)^{n\lambda} \right.\\[.3pc]
\hskip -4pc&\quad\ \left. \times \left(\frac{|\Omega(y-z)| |R_{m+1}(A;x,z)\Vert
f(z)|}{|y-z|^{n-1}|x-z|^m}\right)^2 \chi_{\Gamma(z)}(y,t)\frac{{\rm
d}y{\rm d}t}{t^{n+3}}\right]^{1/2}{\rm d}z \\[.3pc]
\hskip -4pc&\le  C\int_{R^n}\frac{|R_{m+1}(A;x,z)\Vert f(z)|}{|x-z|^m}\\[.3pc]
\hskip -4pc&\quad\ \times \left[\int_0^\infty\int_{|x-y|\le
t}\left(\frac{t}{t+|x-y|}\right) ^{n\lambda}
\frac{\chi_{\Gamma(z)}(y,t)} {(|x-z|-3t)^{2n-2}}\frac{{\rm d}y{\rm
d}t}{t^{n+3}}\right]^{1/2}{\rm d}z \\[.3pc]
\hskip -4pc&\quad\ +C\int_{R^n}\frac{|R_{m+1}(A;x,z)\Vert f(z)|}{|x-z|^m}\\[.3pc]
\hskip -4pc&\quad\ \times \left[\int_0^\infty\sum_{k=0}^\infty\int_{2^kt<|x-y|\le
2^{k+1}t}\left(\frac{t}{t\!+\!|x\!-
\!y|}\right)^{n\lambda}\frac{\chi_{\Gamma(z)}(y,t)t^{-n-3}{\rm d}y{\rm
d}t} {(|x-z|-2^{k+3}t)^{2n-2}}\right]^{1/2}\!{\rm d}z \\[.3pc]
\hskip -4pc&\le C\int_{R_n}\frac{|R_{m+1}(A;x,z)\Vert f(z)|}{|x-z|^{m+1/2}}
\left[\int_{|x-z|/2}^\infty\frac{{\rm
d}t}{(|x-z|-3t)^{2n}}\right]^{1/2}{\rm d}z\\[.3pc]
\hskip -4pc&\quad\ +C\int_{R_n}\frac{|R_{m+1}(A;x,z)\Vert
f(z)|}{|x-z|^{m+1/2}}\\[.3pc]
\hskip -4pc&\quad\ \times \left[\sum_{k=0}^\infty\int_{2^{-2-k}|x-z|}^\infty
2^{-kn\lambda}(2^kt)^nt^{-n} \frac{2^k{\rm
d}t}{(|x-z|-2^{k+3}t)^{2n}}\right]^{1/2}{\rm d}z \\[.3pc]
\hskip -4pc&\le C\int_{R_n}\frac{|R_{m+1}(A;x,z)\Vert f(z)|}{|x-z|^{m+n}}{\rm
d}z\\[.3pc]
\hskip -4pc&\quad\ +C\int_{R_n}\frac{|R_{m+1}(A;x,z)\Vert f(z)|}{|x-z|^{m+n}}{\rm
d}z \left[\sum_{k=0}^\infty 2^{kn(1-\lambda)}\right]^{1/2} \\[.3pc]
\hskip -4pc&= C\int_{R_n}\frac{|R_{m+1}(A;x,z)|}{|x-z|^{m+n}}|f(z)|{\rm d}z.
\end{align*}
Thus, the lemma follows from \cite{2}.
\end{proof}

Now we can prove Theorems~1 and 2. Since $g_\lambda$ and $\mu_\lambda$
are all bounded on $L^p(R^n)$ for $1<p<\infty$ by Lemma~5, it suffices
to verify that $g_\lambda^A$ and $\mu_\lambda^A$ satisfy the size
condition in the Main Theorem.

For $g_\lambda^A$, we write
\begin{align*}
F_t^A(f)(x,y) &= \int \frac{R_{m+1}(\tilde
A;x,z)}{|x-z|^m}\psi_t(y-z)f(z){\rm d}z\\[.2pc]
&= \int \frac{R_{m+1}(\tilde A;x,z)}{|x-z|^m}\psi_t(y-z)f_2(z){\rm d}z\\[.2pc]
&\quad\ +\int\frac{R_m(\tilde A;x,z)}{|x-z|^m}\psi_t(y-z)f_1(z){\rm d}z\\[.2pc]
&\quad\
-\sum_{|\alpha|=m}\frac{1}{\alpha!}\int\frac{(x-z)^\alpha\psi_t(y-z)}{|x-
z|^m}D^\alpha\tilde A(z)f_1(z){\rm d}z,
\end{align*}
then
\begin{align*}
&\frac{1}{|Q|^{1+\beta/n}}\int_Q\left|g_\lambda^A(f)(x)-
g_\lambda^{\tilde A}(f_2)(x_0)\right|{\rm d}x \\[.2pc]
&=\frac{1}{|Q|^{1+\beta/n}}\int_Q\left\vert\left\Vert
\left(\frac{t}{t+|x-y|} \right)^{n\lambda/2}F_t^A(f)(x,y)\right\Vert \right.\\[.2pc]
&\quad\ - \left. \left\Vert
\left(\frac{t}{t+|x-y|} \right)^{n\lambda/2}F_t^{\tilde A}(f)(x_0,y) \right\Vert
\right\vert{\rm d}x \\[.2pc]
&\le \frac{1}{|Q|^{1+\beta/n}}\int_Q
\left\Vert \left(\frac{t}{t+|x-y|} \right)^{n\lambda/2}F_t\left(\frac{R_m(\tilde
A;x,\cdot)}{|x-\cdot|^m}f_1\right)(y)\right\Vert{\rm d}x \\[.2pc]
&\quad\ + \frac{1}{|Q|^{1+\beta/n}}\int_Q\sum_{|\alpha|=
m}\frac{1}{\alpha!}\left\Vert \left(\frac{t}{t+|x-
y|} \right)^{n\lambda/2} \right.\\[.2pc]
&\quad\ \left. \times F_t \left(\frac{(x-\cdot)^{\alpha}}{|x-
\cdot|^m}D^\alpha\tilde A f_1\right)(y) \begin{array}{@{}c@{}}\\ \\[.7pc] \end{array}\right\Vert{\rm d}x \\[.2pc]
&\quad\ +\frac{1}{|Q|^{1+\beta/n}}\int_Q \left\Vert
\left(\frac{t}{t+|x-y|} \right)^{n\lambda/2}F_t^{\tilde
A}(f_2)(x,y) \right.\\[.2pc]
&\quad\ \left. - \left(\frac{t}{t+|x-y|} \right)^{n\lambda/2}F_t^{\tilde
A}(f_2)(x_0, y) \right\Vert {\rm d}x \\[.2pc]
&:=  {\rm I} + {\rm II} + {\rm III}.
\end{align*}
For I and II, similar to the proof of Lemma~5 and the Main Theorem, we get
\begin{align*}
{\rm I} &\le C\sum_{|\alpha|=m}\Vert D^\alpha A\Vert_{\dot\land_\beta}\frac{1}{|Q|}\int_Q|g_\lambda(f_1)(x)|{\rm d}x\\[.2pc]
&\le C\sum_{|\alpha|=m}\Vert D^\alpha A\Vert_{\dot\land_\beta}\Vert g_\lambda(f_1)\Vert_{L^r}|Q|^{-1/r}
\end{align*}
\begin{align*}
&\le C\sum_{|\alpha|=m}\Vert D^\alpha A\Vert_{\dot\land_\beta}\Vert f_1\Vert_{L^r}|Q|^{-1/r}\\[.2pc]
&\le C\sum_{|\alpha|=m}\Vert D^\alpha A\Vert_{\dot\land_\beta}M_r(f)(\tilde x)
\end{align*}
and
\begin{align*}
{\rm II} &\le \frac{C}{|Q|^{1+\beta/n}}\sum_{|\alpha|=m}\Vert
g_\lambda((D^\alpha A-(D^\alpha A)_{\tilde Q})f\chi_{\tilde
Q})\Vert_{L^r}|Q|^{1-1/r} \\[.2pc]
&\le C|Q|^{-\beta/n-1/r}\sum_{|\alpha|=m}\Vert (D^\alpha A-(D^\alpha
A)_{\tilde Q})f\chi_{\tilde Q}\Vert_{L^r} \\[.2pc]
&\le C\sum_{|\alpha|=m}\Vert D^\alpha
A\Vert_{\dot\land_\beta}M_r(f)(\tilde x).
\end{align*}
For III, we write
\begin{align*}
&\left(\frac{t}{t+|x-y|}\right)^{n\lambda/2}F_t^{\tilde A}(f_2)(x,y)-\left(\frac{t}{t+|x_0-y|}\right)^{n\lambda/2}F_t^{\tilde A}(f_2)(x_0,y)  \\[.2pc]
 &= \int_{R^n}\left(\frac{t}{t+|x-y|}\right)^{n\lambda/2}\left[\frac{1}{|x\!-\!z|^m}-\frac{1}{|x_0-z|^m}\right]\\[.2pc]
&\quad\ \times R_m(\tilde A;x,z)\psi_t(y-z)f_2(z){\rm d}z\\[.2pc]
 &\quad\ +\int_{R^n}\left(\frac{t}{t\!+\!|x\!-\!y|}\right)^{n\lambda/2}\frac{\psi_t(y-z)f_2(z)}{|x_0-z|^m}[R_m(\tilde A;x,z)\!-\!R_m(\tilde A;x_0, z)]{\rm d}z\\[.2pc]
 &\quad\ +\int_{R^n}\left[\left(\frac{t}{t+|x-y|}\right)^{n\lambda/2}-\left(\frac{t}{t+|x_0-y|}\right)^{n\lambda/2}\right]\\[.2pc]
&\quad\ \times \frac{R_m(\tilde A;x_0,z)\psi_t(y-z)f_2(z)}{|x_0-z|^m}{\rm d}z\\[.2pc]
 &\quad\ -\sum_{|\alpha|=m}\frac{1}{\alpha!}\int_{R^n}\left[\left(\frac{t}{t+|x-y|}\right)^{n\lambda/2}\frac{(x-z)^\alpha}{|x-z|^m} \right.\\[.2pc]
&\quad\ \left. -\left(\frac{t}{t+|x_0-y|}\right)^{n\lambda/2}\frac{(x_0-z)^\alpha}{|x_0-z|^m}\right]D^\alpha\tilde A(z)\psi_t(y-z)f_2(z){\rm d}z\\[.2pc]
 &:= {\rm III}_1 + {\rm III}_2 + {\rm III}_3 + {\rm III}_4.
\end{align*}
Note that $|x-z|\sim |x_0-z|$ for $x\in Q$ and $z\in R^n\setminus \tilde
Q$. By the condition of $\psi$ and similar to the proof of Lemma~5, we
obtain
\begin{align*}
 &\frac{1}{|Q|^{1+\beta/n}}\int_Q\Vert {\rm III}_1\Vert {\rm d}x\\[.2pc]
 &\le \frac{C}{|Q|^{1+\beta/n}}\int_Q\left(\int_{R^n\setminus \tilde Q}\frac{|x-x_0|}{|x_0-z|^{m+n+1}}
  |R_m(\tilde A;x,z)\Vert f(z)|{\rm d}z\right){\rm d}x
\end{align*}
\begin{align*}
 &\le C\sum_{|\alpha|=m}\Vert D^\alpha A\Vert_{\dot\land_\beta}\sum_{k=0}^\infty
 \int_{2^{k+1}\tilde Q\setminus2^{k+1}\tilde Q}\frac{|x-x_0|}{|x_0-z|^{n+1}}|f(z)|{\rm d}z   \\[.2pc]
 &\le C\sum_{|\alpha|=m}\Vert D^\alpha A\Vert_{\dot\land_\beta}\sum_{k=1}^\infty 2^{-k}\frac{1}{|2^k\tilde Q|}
     \int_{2^k\tilde Q}|f(z)|{\rm d}z   \\[.2pc]
 &\le C\sum_{|\alpha|=m}\Vert D^\alpha A\Vert_{\dot\land_\beta}\sum_{k=1}^\infty 2^{-k}M(f)(\tilde x)   \\[.2pc]
 &\le C\sum_{|\alpha|=m}\Vert D^\alpha A\Vert_{\dot\land_\beta}M(f)(\tilde x).
\end{align*}
 For ${\rm III}_2$, by the formula \cite{5}
\begin{equation*}
R_m(\tilde A; x, z)-R_m(\tilde A; x_0,
z)=\sum_{|\eta|<m}\frac{1}{\eta!}R_{m-|\eta|}(D^\eta\tilde A; x,
x_0)(x-z)^\eta
\end{equation*}
and Lemma~4, we get
\begin{align*}
&|R_m(\tilde A; x, z)-R_m(\tilde A; x_0, z)|\\[.2pc]
&\quad\ \le C\sum_{|\alpha|=m}\Vert
D^\alpha A\Vert_{\dot\land_\beta} |Q|^{\beta/n}|x-x_0\Vert x_0-z|^{m-1}.
\end{align*}
Thus
\begin{align*}
 &\frac{1}{|Q|^{1+\beta/n}}\int_Q\Vert {\rm III}_2\Vert {\rm d}x\\[.3pc]
 &\le C\frac{1}{|Q|^{1+\beta/n}}\int_Q\int_{R^n\setminus\tilde Q}
 \frac{|R_m(\tilde A; x,z)-R_m(\tilde A; x_0, z)|}{|x_0-z|^{m+n}}|f(z)|{\rm d}z {\rm d}x   \\[.3pc]
 &\le C\sum_{|\alpha|=m}\Vert D^\alpha A\Vert_{\dot\land_\beta}\sum_{k=0}^\infty\int_{2^{k+1}\tilde Q\setminus2^k\tilde Q}
 \frac{|x-x_0|}{|x_0-y|^{n+1}}|f(z)|{\rm d}z    \\[.3pc]
 &\le C\sum_{|\alpha|=m}\Vert D^\alpha A\Vert_{\dot\land_\beta}|Q|^{\beta/n}M(f)(\tilde x).
 \end{align*}
For ${\rm III}_3$, by the inequality: $a^{1/2}-b^{1/2}\le (a-b)^{1/2}$
for $a\ge b>0$, we gain
\begin{align*}
\hskip -4.1pc &\frac{1}{|Q|^{1+\beta/n}}\int_Q\Vert {\rm III}_3\Vert {\rm d}x
\le \frac{C}{|Q|^{1+\beta/n}}\\[.2pc]
\hskip -4.1pc &\quad\ \times \int_Q\int_{R^n} \left(\!\int_{R_+^{n+1}}
\left[\frac{t^{n\lambda/2}|x\!-\!x_0|^{1/2}|\psi_t(y\!-\!
z)\vert\vert R_m(\tilde A; x_0, z)\vert\vert
f_2(z)|}{|x\!-\!z|^m(t+|x\!-\!y|)^{(n\lambda+1)/2}}\right]^2
\frac{{\rm d}y {\rm d}t}{t^{n+1}} \!\right)^{1/2}\!{\rm d}z {\rm d}x
\end{align*}\vspace{2pc}
\begin{align*}
\hskip -4.1pc &\le \frac{C}{|Q|^{1+\beta/n}}\int_Q\int_{R^n}\frac{|f_2(z)\vert\vert
R_m(\tilde A; x_0, z)\vert\vert x-x_0|^{1/2}}{|x-z|^m}\\[.2pc]
\hskip -4.1pc &\quad\ \times \left(\int_0^\infty\frac{{\rm
d}t}{(t+|x-z|)^{2n+2}}\right)^{1/2}{\rm d}z {\rm d}x \\[.2pc]
\hskip -4.1pc &\le C\sum_{|\alpha|=m}\Vert D^\alpha
A\Vert_{\dot\land_\beta}\int_{R^n}\frac{|x-x_0|^{1/2}}{|x_0-
z|^{n+1/2}}|f_2(z)|{\rm d}z\\[.2pc]
\hskip -4.1pc &\le C\sum_{|\alpha|=m}\Vert D^\alpha
A\Vert_{\dot\land_\beta}M(f)(\tilde x).
\end{align*}
For ${\rm III}_4$, by Lemma~4, we know that $|D^\alpha A(z)-(D^\alpha
A)_{\tilde Q}|\le \Vert D^\alpha A\Vert_{\dot\land_\beta}|x_0-z|^\beta$.
Thus, similar to the proof of ${\rm III}_1$ and ${\rm III}_3$, we obtain
\begin{align*}
&\frac{1}{|Q|^{1+\beta/n}}\int_Q\Vert {\rm III}_4\Vert {\rm d}x   \\[.2pc]
&\le \frac{C}{|Q|^{1+\beta/n}}\int_Q\sum_{|\alpha|=m}\int_{R^n}\left(\frac{|x-x_0|}{|x_0-z|^{n+1}}
+\frac{|x-x_0|^{1/2}}{|x_0-z|^{n+1/2}}\right)\\[.2pc]
&\quad\ \times |f_2(z)\vert\vert D^\alpha\tilde A(z)|{\rm d}z {\rm d}x   \\[.2pc]
&\le C\sum_{|\alpha|=m}\Vert D^\alpha A\Vert_{\dot\land_\beta}\sum_{k=1}^\infty (2^{k(\beta-1)}+2^{k(\beta-1/2)})M(f)(\tilde x)   \\[.2pc]
&\le C\sum_{|\alpha|=m}\Vert D^\alpha A\Vert_{\dot\land_\beta}M(f)(\tilde x).
\end{align*}
Thus,
\begin{equation*}
{\rm III}\le C\sum_{|\alpha|=m}\Vert D^\alpha
A\Vert_{\dot\land_\beta}M(f)(\tilde x).
\end{equation*}
For $\mu_\lambda^A$, similarly, we have
\begin{align*}
&\frac{1}{|Q|^{1+\beta/n}}\int_Q |\mu_\lambda^A(f)(x)-
\mu_\lambda^{\tilde A}(f_2)(x_0) |{\rm d}x\\[.2pc]
&= \frac{1}{|Q|^{1+\beta/n}}\int_Q \left\vert \left\Vert
\left(\frac{t}{t+|x-y|}\right)^{n\lambda/2}F_t^A(f)(x,y)\right\Vert \right.\\[.2pc]
&\quad\ - \left. \left\Vert \left(\frac{t}{t+|x-y|}
\right)^{n\lambda/2}F_t^{\tilde A}(f_2)(x_0,y)\right\Vert
\right\vert{\rm d}x\\[.2pc]
&\le \frac{1}{|Q|^{1+\beta/n}}\int_Q\left\Vert
\left(\frac{t}{t+|x-y|}\right)^{n\lambda/2} F_t\left(\frac{R_m(\tilde A;
x,\cdot)}{|x-\cdot|^m}f_1\right)(y)\right\Vert {\rm d}x \\[.2pc]
&\quad\ +
\frac{1}{|Q|^{1+\beta/n}}\int_Q\sum_{|\alpha|=m}\frac{1}{\alpha!}
\left\Vert \left(\frac{t}{t\!+\!|x\!-\!y|} \!\right)^{n\lambda/2}\! F_t \!\left(\frac{(x\!-\!
\cdot)^\alpha}{|x\!-\!\cdot|^m} D^\alpha\tilde A
f_1\!\right)(y)\right\Vert{\rm d}x 
\end{align*}
\begin{align*}
&\quad\ +\frac{1}{|Q|^{1+\beta/n}}\int_Q \left\Vert
\left(\frac{t}{t+|x-y|} \right)^{n\lambda/2} F_t^{\tilde
A}(f_2)(x,y) \right.\\[.2pc]
&\quad\ - \left. \left(\frac{t}{t+|x-y|}\right)^{n\lambda/2}F_t^{\tilde
A}(f_2)(x_0,y)\right\Vert {\rm d}x\\[.2pc]
&:= {\rm I} + {\rm II} + {\rm III}
\end{align*}
and
\begin{align*}
{\rm I} &\le C\sum_{|\alpha|=m}\Vert D^\alpha A\Vert_{\dot\land_\beta}\frac{1}{|Q|}\int_Q|\mu_\lambda(f_1)(x)|{\rm d}x\\[.2pc]
&\le C\sum_{|\alpha|=m}\Vert D^\alpha A\Vert_{\dot\land_\beta}\Vert \mu_\lambda(f_1)\Vert_{L^r}|Q|^{-1/r}  \\[.2pc]
&\le C\sum_{|\alpha|=m}\Vert D^\alpha A\Vert_{\dot\land_\beta}\Vert f_1\Vert_{L^r}|Q|^{-1/r}\\[.2pc]
&\le C\sum_{|\alpha|=m}\Vert D^\alpha A\Vert_{\dot\land_\beta}M_r(f)(\tilde x),\\[.2pc]
{\rm II} &\le \frac{C}{|Q|^{1+\beta/n}}\sum_{|\alpha|=m}\Vert \mu_\lambda((D^\alpha A-(D^\alpha A)_{\tilde Q})f\chi_{\tilde Q})\Vert_{L^r}|Q|^{1-1/r}  \\[.2pc]
&\le C|Q|^{-\beta/n-1/r}\sum_{|\alpha|=m}\Vert (D^\alpha A-(D^\alpha A)_{\tilde Q})f\chi_{\tilde Q}\Vert_{L^r}   \\[.2pc]
&\le C\sum_{|\alpha|=m}\Vert D^\alpha A\Vert_{\dot\land_\beta}M_r(f)(\tilde x).
\end{align*}
For III, we write
\begin{align*}
&\left(\frac{t}{t+|x-y|}\right)^{n\lambda/2}F_t^{\tilde A}(f_2)(x,y)-
\left(\frac{t}{t+|x_0-y|}\right)^{n\lambda/2}F_t^{\tilde
A}(f_2)(x_0,y)\\[.2pc]
&= \int_{|y-z|\le
t}\left(\frac{t}{t+|x-y|}\right)^{n\lambda/2}\left[\frac{1}{|x-z|^m}-
\frac{1}{|x_0-z|^m}\right]\\[.2pc]
&\quad\ \times \frac{\Omega(y-z)R_m(\tilde
A;x,z)f_2(z)}{|y-z|^{n-1}}{\rm d}z \\[.2pc]
&\quad\ +\int_{|y-z|\le
t}\left(\frac{t}{t+|x-y|}\right)^{n\lambda/2}\frac{\Omega(y-z)f_2(z)}
{|y-z|^{n-1}|x_0-z|^m}\\[.2pc]
&\quad\ \times [R_m(\tilde A;x,z)-R_m(\tilde A;x_0, z)]{\rm d}z\\[.2pc]
&\quad\ +\int_{|y-z|\le
t}\left[\left(\frac{t}{t+|x-y|}\right)^{n\lambda/2}-
\left(\frac{t}{t+|x_0-y|}\right)^{n\lambda/2}\right]\\[.2pc]
&\quad\ \times \frac{\Omega(y-z)R_m(\tilde
A;x_0,z)f_2(z)}{|y-z|^{n-1}|x_0-z|^m}{\rm d}z
\end{align*}
\begin{align*}
&\quad\ -\sum_{|\alpha|=m}\frac{1}{\alpha!}\int_{|y-z|\le t}\left[\left(\frac{t}{t+|x-y|}\right)^{n\lambda/2}\frac{(x-z)^\alpha}{|x-z|^m} \right.\\[.2pc]
&\quad\ \left. -\left(\frac{t}{t+|x_0-y|}\right)^{n\lambda/2}\frac{(x_0-z)^\alpha}{|x_0-z|^m}\right]
\frac{\Omega(y-z)D^\alpha\tilde A(z)f_2(z)}{|y-z|^{n-1}}{\rm d}z\\[.2pc]
&:= {\rm III}_1 + {\rm III}_2 + {\rm III}_3 + {\rm III}_4.
\end{align*}
By the condition of $\Omega$ and similar to the proof of Lemma~5, we get
\begin{align*}
&\frac{1}{|Q|^{1+\beta/n}}\int_Q\Vert {\rm III}_1\Vert {\rm d}x\\[.2pc]
&\le \frac{C}{|Q|^{1+\beta/n}}\int_Q\left(\int_{R^n\setminus \tilde
Q}\frac{|x-x_0|}{|x_0-z|^{m+n+1}} |R_m(\tilde A;x,z)\Vert
f(z)|{\rm d}z\right){\rm d}x \\[.2pc]
&\le C\sum_{|\alpha|=m}\Vert D^\alpha
A\Vert_{\dot\land_\beta}\sum_{k=0}^\infty \int_{2^{k+1}\tilde
Q\setminus2^{k+1}\tilde Q}\frac{|x-x_0|}{|x_0-z|^{n+1}}|f(z)|{\rm d}z\\[.2pc]
&\le C\sum_{|\alpha|=m}\Vert D^\alpha
A\Vert_{\dot\land_\beta}M(f)(\tilde x),\\[.2pc]
&\frac{1}{|Q|^{1+\beta/n}}\int_Q\Vert {\rm III}_2\Vert {\rm d}x\\[.2pc]
&\le C\frac{1}{|Q|^{1+\beta/n}}\int_Q\int_{R^n\setminus\tilde Q}
\frac{|R_m(\tilde A; x,z)-R_m(\tilde A; x_0, z)|}{|x_0-z|^{m+n}}|f(z)|{\rm d}z
{\rm d}x \\[.2pc]
&\le C\sum_{|\alpha|=m}\Vert D^\alpha
A\Vert_{\dot\land_\beta}\sum_{k=0}^\infty\int_{2^{k+1}\tilde
Q\setminus2^k\tilde Q} \frac{|x-x_0|}{|x_0-y|^{n+1}}|f(z)|{\rm d}z \\[.2pc]
&\le C\sum_{|\alpha|=m}\Vert D^\alpha
A\Vert_{\dot\land_\beta}|Q|^{\beta/n}M(f)(\tilde x),\\[.2pc]
&\frac{1}{|Q|^{1+\beta/n}}\int_Q\Vert {\rm III}_3 \Vert {\rm d}x\\[.2pc]
&\le \frac{C}{|Q|^{1+\beta/n}}\int_Q\int_{R^n}
\left(\int_{R_+^{n+1}}\left[\frac{t^{n\lambda/2}|x-
x_0|^{1/2}\chi_{\Gamma(z)}(y,t)|f_2(z)|}
{(t+|x-y|)^{(n\lambda+1)/2}|y-z|^{n-1}} \begin{array}{c}\\ \\ \\[-.7pc]\end{array}\right. \begin{array}{c}\\ \\ \\[-.1pc] \end{array}\right.\\[.2pc]
&\quad\ \left.\left. \times \frac{|R_m(\tilde
A;x_0,z)|}{|x_0-z|^m}\right]^2\frac{{\rm d}y{\rm
d}t}{t^{n+3}}\right)^{1/2}{\rm d}z{\rm d}x \\[.2pc]
&\le \frac{C}{|Q|^{1+\beta/n}}\int_Q\int_{R^n}\frac{|R_m(\tilde A; x_0,
z)\vert\vert f_2(z)\vert\vert x-x_0|^{1/2}}{|x_0-z|^{m+n+1/2}}{\rm
d}z{\rm d}x \\[.2pc]
&\le C\sum_{|\alpha|=m}\Vert D^\alpha A\Vert_{\dot\land_\beta}M(f)(\tilde x)\\
\end{align*}
and
\begin{align*}
&\frac{1}{|Q|^{1+\beta/n}}\int_Q\Vert {\rm III}_4\Vert {\rm d}x\\[.2pc]
&\le\frac{C}{|Q|^{1+\beta/n}}\int_Q\sum_{|\alpha|=
m}\int_{R^n}\left(\frac{|x-x_0|}{|x_0-z|^{n+1}}
+\frac{|x-x_0|^{1/2}}{|x_0-z|^{n+1/2}}\right)\\[.2pc]
&\quad\ \times |f_2(z)\Vert D^\alpha\tilde A(z)|{\rm d}z {\rm d}x \\[.2pc]
&\le C\sum_{|\alpha|=m}\Vert D^\alpha
A\Vert_{\dot\land_\beta}\sum_{k=1}^\infty
(2^{k(\beta-1)}+2^{k(\beta-1/2)})M(f)(\tilde x) \\[.2pc]
&\le C\sum_{|\alpha|=m}\Vert D^\alpha
A\Vert_{\dot\land_\beta}M(f)(\tilde x).
\end{align*}
Thus,
\begin{equation*}
{\rm III}\le C\sum_{|\alpha|=m}\Vert D^\alpha A\Vert_{\dot\land_\beta}
M(f)(\tilde x).
\end{equation*}
These yield the desired results.

\end{document}